\documentclass[11pt]{article}

\usepackage{graphicx}

\usepackage[utf8x]{inputenc}
\usepackage{tikz}

\usepackage{latexsym}
\usepackage{amssymb}

\newcommand{\excise}[1]{}%{$\star$\textsc{#1}$\star$}

\newtheorem{thm}{Theorem}%[section]

\newtheorem{Example}[thm]{Example}
\newtheorem{Remark}[thm]{Remark}
\newtheorem{Alg}[thm]{Algorithm}
\newtheorem{Defn}[thm]{Definition}

\def\hln{\\[-.2ex]\hline}

        {\begin{list}
                {\noindent\makebox[0mm][r]{\arabic{enumi}.}}
                {\leftmargin=5.5ex \usecounter{enumi}}
        }
        {\end{list}}

        {\begin{list}
                {\noindent\makebox[0mm][r]{(\roman{enumi})}}
                {\leftmargin=5.5ex \usecounter{enumi}}
        }
        {\end{list}}

\def\<{\langle}
\def\>{\rangle}
\def\0{{\mathbf 0}}
\def\1{{\mathbf 1}}

\def\start{{\rm start}}

\font\co=lcircle10
\def\petit#1{{\scriptstyle #1}}

\def\jr{\smash{\raise2pt\hbox{\co \rlap{\rlap{\char'005} \char'007}}
               \raise6pt\hbox{\rlap{\vrule height6.5pt}}
               \raise2pt\hbox{\rlap{\hskip4pt \vrule height0.4pt depth0pt
                width7.7pt}}}}
\def\je{\smash{\raise2pt\hbox{\co \rlap{\rlap{\char'005}
                \phantom{\char'007}}}\raise6pt\hbox{\rlap{\vrule height6pt}}}}
\def\+{\smash{\lower2pt\hbox{\rlap{\vrule height14pt}}
                \raise2pt\hbox{\rlap{\hskip-3pt \vrule height.4pt depth0pt
                width14.7pt}}}}
\def\textcross{\ \smash{\lower4pt\hbox{\rlap{\hskip4.15pt\vrule height14pt}}
                \raise2.8pt\hbox{\rlap{\hskip-3pt \vrule height.4pt depth0pt
                width14.7pt}}}\hskip12.7pt}
\def\textelbow{\ \hskip.1pt\smash{\raise2.8pt%
                \hbox{\co \hskip 4.15pt\rlap{\rlap{\char'005} \char'007}
                \lower6.8pt\rlap{\vrule height3.5pt}
                \raise3.6pt\rlap{\vrule height3.5pt}}
                \raise2.8pt\hbox{%
                  \rlap{\hskip-7.15pt \vrule height.4pt depth0pt width3.5pt}%
                  \rlap{\hskip4.05pt \vrule height.4pt depth0pt width3.5pt}}}
                \hskip8.7pt}
\def\plel{$\begin{tinyrc}{
  \begin{array}{@{}|@{\,}c@{\,}|@{}}
      \hline \scriptscriptstyle  +
    \\\hline \cdot
    \\\hline
  \end{array}
  }\end{tinyrc}$\ }
\def\elel{$\begin{tinyrc}{
  \begin{array}{@{}|@{\,}c@{\,}|@{}}
      \hline \,\cdot\,
    \\\hline \cdot
    \\\hline
  \end{array}
  }\end{tinyrc}$\ }
\def\elpl{$\begin{tinyrc}{
  \begin{array}{@{}|@{\,}c@{\,}|@{}}
      \hline \cdot
    \\\hline \scriptscriptstyle   +
    \\\hline
  \end{array}
  }\end{tinyrc}$\ }
\def\plpl{$\begin{tinyrc}{
  \begin{array}{@{}|@{\,}c@{\,}|@{}}
      \hline \scriptscriptstyle   +
    \\\hline \scriptscriptstyle   +
    \\\hline
  \end{array}
  }\end{tinyrc}$\ }

\usepackage[utf8x]{inputenc}

\usepackage{amssymb}

\title{Mitosis algorithm for Grothendieck polynomials}
\author{D.N. Tyurin}
\begin{document}

\maketitle

\begin{abstract} In this article we will introduce the way to extend the mitosis algorithm for Schubert 
polynomials, developed by Ezra Miller, to the case of Grothendieck polynomials.

\end{abstract}

\section*{Introduction}

$$$$

The main subject of the present note is the extension of the results 
obtained by Ezra Miller in $[3]$ to the case of Grothendieck polynomials.
The mitosis algorithm was firstly introduced in $[5]$ as a combinatorial
rule that allows to compute the coefficients of Schubert polynomials 
inductively in terms of special combinatorial objects called $\emph{pipe dreams}$,
developed by Fomin and Kirillov.

Originally, Schubert polynomials are defined by downward induction on $\emph{weak
Bruhat or-}$ $\emph{der}$ where the induction step is represented by applying of the corresponding
$\emph{divided differ-}$ $\emph{ence operator}$. At the same time, according to the  formula of
Billey, Jockusch, and Stanley, the coefficients of each Schubert polynomial might be 
obtained from the set of diagrams in an $n\times n$ grid called $\emph{reduced pipe dreams}$
(or rc-graphs). Thus, mitosis might be considered as an analogue of the 
applying of the divided difference operator: namely, if $w$ and $v$ are permutations
(and $v<w$ in the weak Bruhat order) then by using this algorithm we can obtain the set
of reduced pipe dreams corresponding to the Schubert polynomial $\mathfrak{S}_{v}(x)$ 
(denoted by $\mathcal{RP}(v)$) from the set $\mathcal{RP}(w)$. The note $[3]$ provides a 
short proof of this fact, based on elementary combinatorial properties of reduced pipe
dreams ($[3]$, Theorem $15$).

In turn, the object of our study --- $\emph{Grothendieck polynomials}$ (denoted by
$\mathfrak{G}_{w}(x)$, where $w$ is an arbitrary permutation) --- might be considered
as a generalization of Schubert polynomials. As well as Schubert polynomials they also can 
be defined inductively by using $\emph{isobaric divided difference operators}$ and in terms
of pipe dreams (only now not necessarily reduced ones). In particular, if $\mathfrak{G}_{w}(x)$
and $\mathfrak{S}_{w}(x)$ are the Schubert and Grothendieck polynomials of a permutation $w$
then $\mathfrak{S}_{w}(x)$ is the sum over all monomials of $\mathfrak{G}_{w}(x)$ of the 
minimal degree. Thus, it is reasonable to suggest that mitosis is somehow applied to 
the case of Grothendieck polynomials (i.e. to pipe dreams in general). In the first part of our
note we will describe the way mitosis acts on pipe dreams in general. The second part
will be devoted to the theorem, analogous to the Theorem $15$ of $[5]$. For that we will 
slightly modify the original mitosis algorithm. %(Здесь будет пассаж про наши результаты)%

The plan of the paper is as follows: in the first section we will
give the inductive definitions of Schubert and Grothendieck polynomials, introduce to the 
reader the concept of pipe dreams and describe the way Schubert and Grohendieck polynomials
can be defined combinatorially. The second section is devoted to the original mitosis
algorithm: we give a brief description of this algorithm and prove theorem $2.1$ describing 
the way mitosis acts on the set of pipe dreams $\mathcal{P}(w)$ corresponding to a  
permutation $w$. Section 3, analogously to
$[5]$, provides a special involution on $\mathcal{P}(w)$ that is cruical for the proof
of the main theorem. The final section is devoted to the main theorem of the paper (Theorem $4.1$), which
gives us  inductive way to construct Grothendieck polynomials in terms of the pipe dreams.

 $$$$
 %%%%%%%%%%%%%%%%%%%%%%%%%%%%%%%%%%%%%%%%%%%%%%%%%%%%%%%%%%%%%%%%%%%%%%%%%
\section{Pipe dreams}%%%%%%%%%%%%%%%%%%%%%%%%%%%%%%%%%%%%%%%%%%%%%%%%%%%%
%%%%%%%%%%%%%%%%%%%%%%%%%%%%%%%%%%%%%%%%%%%%%%%%%%%%%%%%%%%%%%%%%%%%%%%%%
$$$$

In this section we will
give the inductive definitions of Schubert and Grothendieck polynomials, and introduce to the 
reader the concept of pipe dreams.

Denote by $s_{i}=(i,i+1)$ the corresponding elementary transposition of $S_{n}$.
It is well known that the set $\{s_{i} | i=1,\ldots,n-1\}$ generates the group $S_{n}$;
with the following relations:

$$s_{i}^{2}=1$$
$$s_{i}s_{j}=s_{j}s_{i} \ \textmd{if} \ |i-j|>1$$
$$s_{i}s_{i+1}s_{i}=s_{i+1}s_{i}s_{i+1}.$$

Thus, we can say that the set $\{s_{i}| i\in\mathbb{N}\}$ generates the group $S_{\infty}$
with the relations described above. Let $w$ be an arbitrary element of $S_{\infty}$.
Then we can define a sequence $a_{1},\ldots, a_{k}$ of minimal length
such that $w=s_{a_{1}}\cdots s_{a_{k}}$. The number $k$ is called the $\emph{length}$ of
$w$ and detoted by $l(w)$  (remark that the sequence itself can be defined in more than one
way).

Let $w$ be an arbitrary element of $S_{\infty}$ and $s_{i}$ be an elementary transposition.
Then we say that $ws_{i}>w$ if $l(ws_{i})=l(w)+1$. Otherwise we have $l(ws_{i})=l(w)-1$
which means that $l(w)=l((ws_{i})s_{i})=l(ws_{i})+1$ and, accordingly, $w>ws_{i}$. Now, by
using the property of transivity we can introduce on $S_{\infty}$ a partial order. It is well
known that for each $ n\in\mathbb{N}$ the greatest element of $S_{n}$ is the
$\emph{order reversing permutation}$, denoted by $w_{0}^{(n)}=(n\ n-1\cdots 1)$

Now we can define Schubert and Grothendieck polynomials inductively. For
that we will also need to introduce the sets of $\emph{divided difference operators}$,
denoted by $\partial_{i}$, and $\emph{isobaric divided difference operators}$,
denoted by $\pi_{i}$. These linear operators act on the ring $\mathbb{Z}[x_{1}, x_{2}, . . .]$ as follows:

$$\forall f \in \mathbb{Z}[x_{1}, x_{2}, . . .], \ \partial_{i}(f)=\frac{f-s_{i}(f)}{x_{i}-x_{i+1}}, \ \pi_{i}(f)=\partial_{i}(f-x_{i+1}f). \eqno(1.1)$$

Here $s_{i}(f)$ is a polynomial, obtained from $f$ by interchanging the variables $x_{i}$ and $x_{i+1}$.
Obviously, the result of applying of each operator is also a polynomial with integer coefficients.

$\textbf{Definition 1.1:}$ For an arbitrary element $w$ of $S_{\infty}$ the corresponding Schubert
polynomial can be defined inductively in compliance with the following rules:

1) For the order reversing permutation $w_{0}^{(n)}$ the following equality holds:
$$\mathfrak{S}_{w_{0}^{(n)}}(x)=x_{1}^{n-1}x_{2}^{n-2}\cdots x_{n-1};$$

2)
$$\partial_{i}(\mathfrak{S}_{w}(x))=\mathfrak{S}_{ws_{i}}(x), \ \textmd{if} \ l(ws_{i})=l(w)-1.$$

The corresponding Grothendieck polynomial is defined practically in the same way. Namely:

1) For $w_{0}^{(n)}$ the Grothendieck polynomial  $\mathfrak{G}_{w_{0}^{(n)}}(x)$
equals the corresponding Schubert polynomial:

$$\mathfrak{G}_{w_{0}^{(n)}}(x)=\mathfrak{S}_{w_{0}^{(n)}}(x)=x_{1}^{n-1}x_{2}^{n-2}\cdots x_{n-1};$$

 2)
 $$\pi_{i}(\mathfrak{G}_{w}(x))=\mathfrak{G}_{ws_{i}}(x), \ \textmd{if} \ l(ws_{i})=l(w)-1.$$

Now we will give the definition of a $\emph{pipe dream}$. Consider the direct product $\mathbb{Z}_{>0}\times\mathbb{Z}_{>0}$,
represented in the form of a table extending infinitely south and east (the box, located in  $i$-th row
and $j$-th column  indexed by pair $(i, j)$). Then a pipe dream is a finite subset of
$\mathbb{Z}_{>0}\times\mathbb{Z}_{>0}$, with it's elements marked with the $+$ symbol
(see the example below).

$\textbf{Example 1.1:}$ The following pipe dream represents the set $\{(1,1);
(1,2);(1,5);(2,2);(3,2);$ $(5,1)\}$.

$$\begin{array}{|c|c|c|c|c|c|}
\multicolumn{6}{c}{}
\\\hline\!+\!&\!+\!&\ \, & \ \,&\!+\!&\ \,
\\\hline      &\!+\!&     &     &       &
\\\hline      &\!+\!&     &     &       &
\\\hline      &      &      &     &       &
\\\hline\!+\!&     &      &     &       &
\\\hline      &      &     &      &       &
\\\hline
\end{array}$$
$$$$

Now in every box with $+$ we put the symbol $\textcross$,
and in every empty box --- the symbol  $\textelbow$. Thus, we obtain the network of strands , crossing
each other at the positions belonging to the pipe dream and avoiding each other at other positions.
The pipe dream is $\emph{reduced}$ if each pair of strands crosses at most once.
Now for every $w\in S_{\infty}$ we define the set $\mathcal{RP}(w)$ of reduced pipe dreams, such
that for each element $D$ of this set the strand entering the $i$-th row, exits
from the $w(i)$-th column.

$\textbf{Example 1.2:}$, The depicted pipe dream is reduced and corresponds to the permutation $w=(261354)$.

$$\begin{array}{ccccccc}
      &\petit2&\petit 6&\petit 1&\petit 3&\petit 5&\petit 4\\
\petit1&  \+ &   \+   &   \jr  &   \jr   &   \+   &   \je\\
\petit2&  \jr &   \+   &   \jr  &   \jr   &   \je   &\\
\petit3&  \jr  &  \+   &   \jr  &   \je   &         &\\
\petit4&  \jr &   \jr   &   \je  &         &          &\\
\petit5&  \+  &  \je  &         &         &          &\\
\petit6&  \je  &       &         &         &          &\\
\end{array}$$

Note that for simplicity we usually do not draw the "sea" of wavy strands, whose entries and exits
are equal.

Now Schubert polynomials can be defined combinatorially. Namely, for an arbitrary pipe dream $D$
we will introduce the following notation

$$x_{D}=\prod\limits_{(i,j)\in D}x_{i}.$$

The following proposition is a nontrivial theorem, which proof will be not be given here
(See $[4]$, Theorem $1.1$).

$\textbf{Theorem 1.1:}$ For an arbitrary element $w$ of  $S_{\infty}$ the following equality holds:

$$\mathfrak{S}_{w}(x)=\sum\limits_{D\in\mathcal{RP}(w)}x_{D}. \eqno(1.2)$$

As a cosequence, the coefficients of the polynomial $\mathfrak{S}_{w}(x)$ are positive.

Consider now an arbitrary pipe dream $B$, whose strands can cross each other more than once.
For each element $(i,j)$ of $D$ we will define its $\emph{antidiagonal index}$ by number
$i+j-1$. Then, by moving across the table from right to left and from top to bottom and
assotiating to each element of $B$ its antidiagonal index, we will obtain a sequence
$(i_{1},\ldots,i_{k})$, where $k$ is the number of crosses of $B$. The corresponding permutation $w$
is produced by multiplying the elementary transpositions $s_{i_{1}}\cdots s_{i_{k}}$ in compliance
with the following rules:

$$s_{i}^{2}=s_{i}$$
$$s_{i}s_{j}=s_{j}s_{i} \ if \ |i-j|>1$$
$$s_{i}s_{i+1}s_{i}=s_{i+1}s_{i}s_{i+1}.$$

In other words, we omit transpositions that decrease length
(the corresponding operation is called $\emph{the Demazure product}$). The set of all pipe dreams, whose Demazure
products are equal to $w$, is denoted by $\mathcal{P}(w)$. The set
$\mathcal{RP}(w)$ is a subset of $\mathcal{P}(w)$. Also in case of a reduced pipe dream the Demazure product
is equivalent to the standard group operation.

$\textbf{Example 1.3:}$ $B$ is a nonreduced pipe dream, belonging to $\mathcal{P}(1423)$.
$$\begin{array}{c}\\\\\\\\\\\raisebox{2ex}{$B \ = \ $}\\\\\\\\\end{array}
\begin{array}{|c|c|c|c|}
\multicolumn{4}{c}{}
\\\hline      &\!+\!&\ \, & \ \,
\\\hline\!+\!&\!+\!&     &
\\\hline\!+\!&      &     &
\\\hline      &      &      &
\\\hline
\end{array}$$
The following theorem is a generalization of Theorem $1.1$ for the case of Grothendeick
polynomials (See $[5]$, Theorem $A$).

$\textbf{Theorem 1.2:}$ For an arbitrary element $w$ of $S_{\infty}$
the following equality holds:
$$\mathfrak{G}_{w}(x)=\sum\limits_{B\in\mathcal{P}(w)}(-1)^{|B|-l(w)}x_{B}. \eqno(1.3)$$

Here by $|B|$ we denote the number of crosses of $B$.

$$$$
 %%%%%%%%%%%%%%%%%%%%%%%%%%%%%%%%%%%%%%%%%%%%%%%%%%%%%%%%%%%%%%%%%%%%%%%%%
\section{Mitosis algorithm}%%%%%%%%%%%%%%%%%%%%%%%%%%%%%%%%%%%%%%%%%%%%%%%%%%%%
%%%%%%%%%%%%%%%%%%%%%%%%%%%%%%%%%%%%%%%%%%%%%%%%%%%%%%%%%%%%%%%%%%%%%%%%%
$$$$
The algorithm of mitosis recursion was developed by Ezra Miller and represents combi-
nation of inductive and combinatorial definitions of Schubert polynomials. More precisely,
it allows us to obtain the set $\mathcal{RP}(ws_{i})$ of reduced pipe-dreams of $ws_{i}$ from
the set $\mathcal{RP}(w)$ if $l(ws_{i})=l(w)-1$. Here we will briefly recall the algorithm.

Consider an arbitrary pipe dream $D$. Introduce the following notiations:

$$\emph{start}_{i}(D)=\min(\{j | (i,j)\notin D\}\cup\{n+1\}).\eqno(2.1)$$

In other words, $\emph{start}_{i}(D)$ is the maximal column index, such that each box in the $i$-th row,
located strictly to the left from $(i,\emph{start}_{i}(D))$ is marked with cross. Also by
$\mathcal{J}_{i}(D)$ we denote the subset of the numbers of columns, located strictly to the left from $(i,\emph{start}_{i}(D))$,
such that for every $j\in\mathcal{J}_{i}(D)$ the box $(i+1,j)$ is empty.

Then each element $p$ of $\mathcal{J}_{i}(D)$ can be associated with a new pipe dream $D_{p}$,constructed
in the following way: first the cross in the box $(i,p)$ is deleted from $D$, then every cross in $i$-th row,
located to the left from $(i,p)$ with in's column index, belonging to $\mathcal{J}_{i}(D)$ is moved down to the
empty box below it.

Now we can define the $\emph{mitosis operator}$.

 $\textbf{Definition 2.1:}$ The $i$-th mitosis operator (denoted by $mitosis_{i}(D)$) sends $D$
to the set $\{D_{p} | p\in\mathcal{J}_{i}(D)\}$.

 $\textbf{Example 2.1:}$
$$
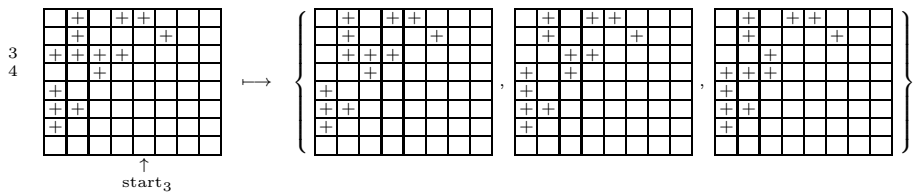
\begin{tinyrc}{
\begin{array}{@{}cc@{}}{}\\\\3&\\4\\\\\\\\\\\\\\\end{array}
\begin{array}{@{}|@{\,}c@{\,}|@{\,}c@{\,}|@{\,}c@{\,}|@{\,}c@{\,}
		 |@{\,}c@{\,}|@{\,}c@{\,}|@{\,}c@{\,}|@{\,}c@{\,}|@{}}
\hline     &  +  &     &  +  &  +  &     &     &
\hln       &  +  &     &     &     &  +  &     &\phantom{+}
\hln    +  &  +  &  +  &  +  &     &     &     &
\hln       &     &  +  &     &     &     &     &
\hln    +  &     &     &     &     &     &     &
\hln    +  &  +  &     &     &     &     &     &
\hln    +  &     &     &     &     &     &     &
\hln       &     &     &     &     &     &\phantom{+}&
\\\hline\multicolumn{4}{@{}c@{}}{}&\multicolumn{1}{@{}c@{}}{\!\!\uparrow}
\\\multicolumn{4}{c}{}&\multicolumn{1}{@{}c@{}}{\makebox[0pt]{$\start_3$}}
\end{array}
\begin{array}{@{\quad}c@{\ }}
  \longmapsto\\\mbox{}\\\mbox{}
\end{array}
\begin{array}{c}
  \left\{\
  \begin{array}{@{}|@{\,}c@{\,}|@{\,}c@{\,}|@{\,}c@{\,}|@{\,}c@{\,}
  		 |@{\,}c@{\,}|@{\,}c@{\,}|@{\,}c@{\,}|@{\,}c@{\,}|@{}}
  \hline     &  +  &     &  +  &  +  &     &     &
  \hln       &  +  &     &     &     &  +  &     &\phantom{+}
  \hln       &  +  &  +  &  +  &     &     &     &
  \hln       &     &  +  &     &     &     &     &
  \hln    +  &     &     &     &     &     &     &
  \hln    +  &  +  &     &     &     &     &     &
  \hln    +  &     &     &     &     &     &     &
  \hln       &     &     &     &     &     &\phantom{+}&
  \\\hline
  \end{array}
  \ ,\
  \begin{array}{@{}|@{\,}c@{\,}|@{\,}c@{\,}|@{\,}c@{\,}|@{\,}c@{\,}
  		 |@{\,}c@{\,}|@{\,}c@{\,}|@{\,}c@{\,}|@{\,}c@{\,}|@{}}
  \hline     &  +  &     &  +  &  +  &     &     &
  \hln       &  +  &     &     &     &  +  &     &\phantom{+}
  \hln       &     &  +  &  +  &     &     &     &
  \hln    +  &     &  +  &     &     &     &     &
  \hln    +  &     &     &     &     &     &     &
  \hln    +  &  +  &     &     &     &     &     &
  \hln    +  &     &     &     &     &     &     &
  \hln       &     &     &     &     &     &\phantom{+}&
  \\\hline
  \end{array}
  \ ,\
  \begin{array}{@{}|@{\,}c@{\,}|@{\,}c@{\,}|@{\,}c@{\,}|@{\,}c@{\,}
  		 |@{\,}c@{\,}|@{\,}c@{\,}|@{\,}c@{\,}|@{\,}c@{\,}|@{}}
  \hline     &  +  &     &  +  &  +  &     &     &
  \hln       &  +  &     &     &     &  +  &     &\phantom{+}
  \hln       &     &  +  &     &     &     &     &
  \hln    +  &  +  &  +  &     &     &     &     &
  \hln    +  &     &     &     &     &     &     &
  \hln    +  &  +  &     &     &     &     &     &
  \hln    +  &     &     &     &     &     &     &
  \hln       &     &     &     &     &     &\phantom{+}&
  \\\hline
  \end{array}
  \ \right\}
\\
  \begin{array}{c}
  \mbox{}\\
  \mbox{}\\
  \end{array}
\end{array}
}\end{tinyrc}
$$

Here $i=3$, and $\mathcal{J}_{i}(D)$ contains 1, 2 and 4.

If $\mathcal{C}$ is a set of pipe dreams, then by $mitosis_{i}(\mathcal{C})$ we mean
the union $\bigcup\limits_{D\in\mathcal{C}}mitosis_{i}(D)$
over all elements of $\mathcal{C}$. If $\mathcal{J}_{i}(D)=\emptyset$, then
 $mitosis_{i}(D)=\emptyset$.

Now consider an element $w$ of $S_{\infty}$ and elemental transposition $s_{i}$,
such that $l(ws_{i})=l(w)-1$.The main result of the corresponding article by E. Miller ($[2]$) is the
following statement:

$\textbf{Theorem 2.1:}$ Disjoint union $\bigsqcup\limits_{D\in\mathcal{RP}(w)}mitosis_{i}(D)$
coincides with the set $\mathcal{RP}(ws_{i})$ of reduced pipe dreams of the permutation $ws_{i}$.

Together with the formula $(1.3)$ this statement gives as the capability to obtain Schubert polynomials
inductively, in terms of pipe dreams. Here we will not give the proof of this theorem (see $[2]$).
Nevertheless, further  we will use the same methods and terms, slightly modifying them.

Let's start from the fact that theorem $\textbf{2.1}$ states, that mitosis operator sends reduced pipe dreams
from $\mathcal{RP}(w)$ to $\mathcal{RP}(ws_{i})$, so it is reasonable to assume, that
the same is right for all elements of $\mathcal{P}(w)$. Nevertheless, it is not quite true.
For our further reasonings we will have to introduce the family of operations of the pipe dreams,
which, in accordance with the original article ($[2]$) will be called $\emph{chute moves}$.

Henceforth, we will distinguish between three types of chute moves. Under the action of a chute move some elements
of a pipe dream change places and some are deleted in a certain way. These elements
are located in two neighbouring rows of the pipe dream.

$\textbf{Definition 2.2:}$ All three types of chute moves (which henceforth will be indexed by numbers 1,2 and 3)
can be defined
graphically in the following way:

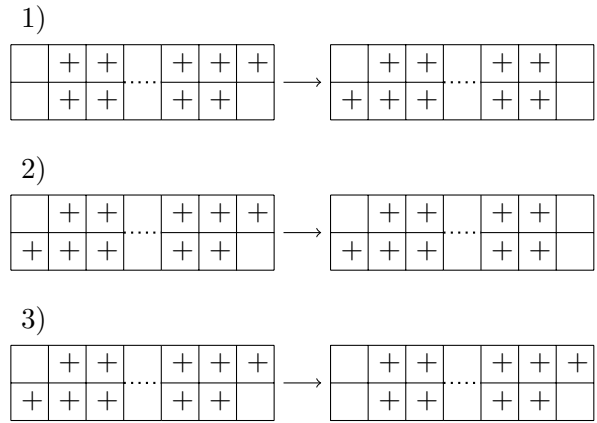
\begin{figure}[!h]
 \centering

\begin{tikzpicture}

\draw[step=0.5cm]  (0,0) grid (1.5,1);
\draw[step=0.5cm]  (2,0) grid (3.5,1);
\draw (2,1)--(2,0);
\draw (1.5,1)--(2,1);
\draw (1.5,0)--(2,0);
\draw[dotted, thick] (1.5,0.5)--(2,0.5);
\draw (0.5,0.5) node [minimum width=0.5cm, minimum height=0.5cm, above right] {$+$};
\draw (0.5,0) node [minimum width=0.5cm, minimum height=0.5cm, above right] {$+$};
\draw (1,0.5) node [minimum width=0.5cm, minimum height=0.5cm, above right] {$+$};
\draw (1,0) node [minimum width=0.5cm, minimum height=0.5cm, above right] {$+$};
\draw (2,0.5) node [minimum width=0.5cm, minimum height=0.5cm, above right] {$+$};
\draw (2,0) node [minimum width=0.5cm, minimum height=0.5cm, above right] {$+$};
\draw (2.5,0.5) node [minimum width=0.5cm, minimum height=0.5cm, above right] {$+$};
\draw (2.5,0) node [minimum width=0.5cm, minimum height=0.5cm, above right] {$+$};
\draw (3,0.5) node [minimum width=0.5cm, minimum height=0.5cm, above right] {$+$};
\draw[xshift=4.25cm, step=0.5cm]  (0,0) grid (1.5,1);
\draw[xshift=4.25cm, step=0.5cm]  (2,0) grid (3.5,1);
\draw[xshift=4.25cm] (2,1)--(2,0);
\draw[xshift=4.25cm] (1.5,1)--(2,1);
\draw[xshift=4.25cm](1.5,0)--(2,0);
\draw[xshift=4.25cm, dotted, thick] (1.5,0.5)--(2,0.5);
\draw[xshift=4.25cm] (0.5,0.5) node [minimum width=0.5cm, minimum height=0.5cm, above right] {$+$};
\draw[xshift=4.25cm] (0.5,0) node [minimum width=0.5cm, minimum height=0.5cm, above right] {$+$};
\draw[xshift=4.25cm] (1,0.5) node [minimum width=0.5cm, minimum height=0.5cm, above right] {$+$};
\draw[xshift=4.25cm] (1,0) node [minimum width=0.5cm, minimum height=0.5cm, above right] {$+$};
\draw[xshift=4.25cm] (2,0.5) node [minimum width=0.5cm, minimum height=0.5cm, above right] {$+$};
\draw[xshift=4.25cm] (2,0) node [minimum width=0.5cm, minimum height=0.5cm, above right] {$+$};
\draw[xshift=4.25cm] (2.5,0.5) node [minimum width=0.5cm, minimum height=0.5cm, above right] {$+$};
\draw[xshift=4.25cm] (2.5,0) node [minimum width=0.5cm, minimum height=0.5cm, above right] {$+$};
\draw[xshift=4.25cm](0,0) node [minimum width=0.5cm, minimum height=0.5cm, above right] {$+$};
\draw[->] (3.625,0.5)--(4.125,0.5);
\draw[yshift=-2cm, ->] (3.625,0.5)--(4.125,0.5);
\draw[yshift=-4cm, ->] (3.625,0.5)--(4.125,0.5);
\draw[yshift=-2cm, step=0.5cm]  (0,0) grid (1.5,1);
\draw[yshift=-2cm, step=0.5cm]  (2,0) grid (3.5,1);
\draw[yshift=-2cm] (2,1)--(2,0);
\draw[yshift=-2cm] (1.5,1)--(2,1);
\draw[yshift=-2cm](1.5,0)--(2,0);
\draw[yshift=-2cm, dotted, thick] (1.5,0.5)--(2,0.5);
\draw[yshift=-2cm] (0.5,0.5) node [minimum width=0.5cm, minimum height=0.5cm, above right] {$+$};
\draw[yshift=-2cm] (0.5,0) node [minimum width=0.5cm, minimum height=0.5cm, above right] {$+$};
\draw[yshift=-2cm] (1,0.5) node [minimum width=0.5cm, minimum height=0.5cm, above right] {$+$};
\draw[yshift=-2cm] (1,0) node [minimum width=0.5cm, minimum height=0.5cm, above right] {$+$};
\draw[yshift=-2cm] (2,0.5) node [minimum width=0.5cm, minimum height=0.5cm, above right] {$+$};
\draw[yshift=-2cm] (2,0) node [minimum width=0.5cm, minimum height=0.5cm, above right] {$+$};
\draw[yshift=-2cm] (2.5,0.5) node [minimum width=0.5cm, minimum height=0.5cm, above right] {$+$};
\draw[yshift=-2cm] (2.5,0) node [minimum width=0.5cm, minimum height=0.5cm, above right] {$+$};
\draw[yshift=-2cm](0,0) node [minimum width=0.5cm, minimum height=0.5cm, above right] {$+$};
\draw[yshift=-2cm](3,0.5) node [minimum width=0.5cm, minimum height=0.5cm, above right] {$+$};

\draw[yshift=-4cm, step=0.5cm]  (0,0) grid (1.5,1);
\draw[yshift=-4cm, step=0.5cm]  (2,0) grid (3.5,1);
\draw[yshift=-4cm] (2,1)--(2,0);
\draw[yshift=-4cm] (1.5,1)--(2,1);
\draw[yshift=-4cm](1.5,0)--(2,0);
\draw[yshift=-4cm, dotted, thick] (1.5,0.5)--(2,0.5);
\draw[yshift=-4cm] (0.5,0.5) node [minimum width=0.5cm, minimum height=0.5cm, above right] {$+$};
\draw[yshift=-4cm] (0.5,0) node [minimum width=0.5cm, minimum height=0.5cm, above right] {$+$};
\draw[yshift=-4cm] (1,0.5) node [minimum width=0.5cm, minimum height=0.5cm, above right] {$+$};
\draw[yshift=-4cm] (1,0) node [minimum width=0.5cm, minimum height=0.5cm, above right] {$+$};
\draw[yshift=-4cm] (2,0.5) node [minimum width=0.5cm, minimum height=0.5cm, above right] {$+$};
\draw[yshift=-4cm] (2,0) node [minimum width=0.5cm, minimum height=0.5cm, above right] {$+$};
\draw[yshift=-4cm] (2.5,0.5) node [minimum width=0.5cm, minimum height=0.5cm, above right] {$+$};
\draw[yshift=-4cm] (2.5,0) node [minimum width=0.5cm, minimum height=0.5cm, above right] {$+$};
\draw[yshift=-4cm](0,0) node [minimum width=0.5cm, minimum height=0.5cm, above right] {$+$};
\draw[yshift=-4cm](3,0.5) node [minimum width=0.5cm, minimum height=0.5cm, above right] {$+$};

\draw[yshift=-4cm, xshift=4.25cm, step=0.5cm]  (0,0) grid (1.5,1);
\draw[yshift=-4cm, xshift=4.25cm, step=0.5cm]  (2,0) grid (3.5,1);
\draw[yshift=-4cm, xshift=4.25cm] (2,1)--(2,0);
\draw[yshift=-4cm, xshift=4.25cm] (1.5,1)--(2,1);
\draw[yshift=-4cm, xshift=4.25cm](1.5,0)--(2,0);
\draw[yshift=-4cm, xshift=4.25cm, dotted, thick] (1.5,0.5)--(2,0.5);
\draw[yshift=-4cm, xshift=4.25cm] (0.5,0.5) node [minimum width=0.5cm, minimum height=0.5cm, above right] {$+$};
\draw[yshift=-4cm, xshift=4.25cm] (0.5,0) node [minimum width=0.5cm, minimum height=0.5cm, above right] {$+$};
\draw[yshift=-4cm, xshift=4.25cm] (1,0.5) node [minimum width=0.5cm, minimum height=0.5cm, above right] {$+$};
\draw[yshift=-4cm, xshift=4.25cm] (1,0) node [minimum width=0.5cm, minimum height=0.5cm, above right] {$+$};
\draw[yshift=-4cm, xshift=4.25cm] (2,0.5) node [minimum width=0.5cm, minimum height=0.5cm, above right] {$+$};
\draw[yshift=-4cm, xshift=4.25cm] (2,0) node [minimum width=0.5cm, minimum height=0.5cm, above right] {$+$};
\draw[yshift=-4cm, xshift=4.25cm] (2.5,0.5) node [minimum width=0.5cm, minimum height=0.5cm, above right] {$+$};
\draw[yshift=-4cm, xshift=4.25cm] (2.5,0) node [minimum width=0.5cm, minimum height=0.5cm, above right] {$+$};
\draw[yshift=-4cm, xshift=4.25cm](3,0.5) node [minimum width=0.5cm, minimum height=0.5cm, above right] {$+$};

\draw[yshift=-2cm, xshift=4.25cm, step=0.5cm]  (0,0) grid (1.5,1);
\draw[yshift=-2cm, xshift=4.25cm, step=0.5cm]  (2,0) grid (3.5,1);
\draw[yshift=-2cm, xshift=4.25cm] (2,1)--(2,0);
\draw[yshift=-2cm, xshift=4.25cm] (1.5,1)--(2,1);
\draw[yshift=-2cm, xshift=4.25cm](1.5,0)--(2,0);
\draw[yshift=-2cm, xshift=4.25cm, dotted, thick] (1.5,0.5)--(2,0.5);
\draw[yshift=-2cm, xshift=4.25cm] (0.5,0.5) node [minimum width=0.5cm, minimum height=0.5cm, above right] {$+$};
\draw[yshift=-2cm, xshift=4.25cm] (0.5,0) node [minimum width=0.5cm, minimum height=0.5cm, above right] {$+$};
\draw[yshift=-2cm, xshift=4.25cm] (1,0.5) node [minimum width=0.5cm, minimum height=0.5cm, above right] {$+$};
\draw[yshift=-2cm, xshift=4.25cm] (1,0) node [minimum width=0.5cm, minimum height=0.5cm, above right] {$+$};
\draw[yshift=-2cm, xshift=4.25cm] (2,0.5) node [minimum width=0.5cm, minimum height=0.5cm, above right] {$+$};
\draw[yshift=-2cm, xshift=4.25cm] (2,0) node [minimum width=0.5cm, minimum height=0.5cm, above right] {$+$};
\draw[yshift=-2cm, xshift=4.25cm] (2.5,0.5) node [minimum width=0.5cm, minimum height=0.5cm, above right] {$+$};
\draw[yshift=-2cm, xshift=4.25cm] (2.5,0) node [minimum width=0.5cm, minimum height=0.5cm, above right] {$+$};
\draw[yshift=-2cm, xshift=4.25cm](0,0) node [minimum width=0.5cm, minimum height=0.5cm, above right] {$+$};

\draw (0,1) node [above right] {$1)$};
\draw (0,-1) node [above right] {$2)$};
\draw (0,-3) node [above right] {$3)$};
\end{tikzpicture}
\end{figure}

As we can see, for any chute move we can uniquely describe an inverse one.
The following statement is right for both chute moves and inverse chute moves:

$\textbf{Theorem 2.2:}$ If a pipe dream $D$ belongs to the set $\mathcal{P}(w)$, then the result
of applying of any of the chute moves 1-3 also belongs to $\mathcal{P}(w)$.

$\textbf{Proof:}$ We will prove the statement only for the cases of chute move -1 and chute move-2,
since chute move-3 can be introduced as a result of consistent applying of chute moves 2 and 1.

1) Let us suppose that  the $m$-th and $(m+1)$-th rows of $D$ look like:

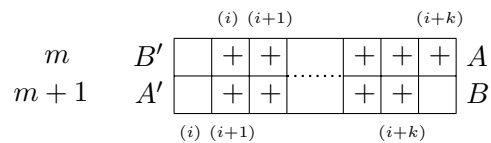
\begin{figure}[!ht]
\centering
\begin{tikzpicture}
\draw[step=0.5cm]  (0,0) grid (1.5,1);
\draw[xshift=2.25cm, step=0.5cm]  (0,0) grid (1.5,1);
\draw (1.5,1)--(2.25,1);
\draw (1.5,0)--(2.25,0);
\draw[dotted, thick] (1.5,0.5)--(2.25,0.5);
\draw (0,0) node [minimum width=0.5cm, minimum height=0.5cm, above left] {$A'$};
\draw (0,0.5) node [minimum width=0.5cm, minimum height=0.5cm, above left] {$B'$};
\draw (3.75,0) node [minimum width=0.5cm, minimum height=0.5cm, above right] {$B$};
\draw (3.75,0.5) node [minimum width=0.5cm, minimum height=0.5cm, above right] {$A$};

\draw (0.5,0) node [minimum width=0.5cm, minimum height=0.5cm, below left] {$\scriptscriptstyle(i)$};
\draw (0.33,0) node [minimum width=0.5cm, minimum height=0.5cm, below right] {$\scriptscriptstyle(i+1)$};

\draw (1,1) node [minimum width=0.5cm, minimum height=0.5cm, above left] {$\scriptscriptstyle(i)$};
\draw (0.83,1) node [minimum width=0.5cm, minimum height=0.5cm, above right] {$\scriptscriptstyle(i+1)$};

\draw (2.58,0) node [minimum width=0.5cm, minimum height=0.5cm, below right] {$\scriptscriptstyle(i+k)$};
\draw (3.08,1) node [minimum width=0.5cm, minimum height=0.5cm, above right] {$\scriptscriptstyle(i+k)$};

\draw (0.5,0.5) node [minimum width=0.5cm, minimum height=0.5cm, above right] {$+$};
\draw (0.5,0) node [minimum width=0.5cm, minimum height=0.5cm, above right] {$+$};
\draw (1,0.5) node [minimum width=0.5cm, minimum height=0.5cm, above right] {$+$};
\draw (1,0) node [minimum width=0.5cm, minimum height=0.5cm, above right] {$+$};

\draw[xshift=1.75cm] (0.5,0.5) node [minimum width=0.5cm, minimum height=0.5cm, above right] {$+$};
\draw[xshift=1.75cm] (0.5,0) node [minimum width=0.5cm, minimum height=0.5cm, above right] {$+$};
\draw[xshift=1.75cm] (1,0.5) node [minimum width=0.5cm, minimum height=0.5cm, above right] {$+$};
\draw[xshift=1.75cm] (1,0) node [minimum width=0.5cm, minimum height=0.5cm, above right] {$+$};
\draw[xshift=2.75cm] (0.5,0.5) node [minimum width=0.5cm, minimum height=0.5cm, above right] {$+$};

\draw (-1,0) node [minimum width=0.5cm, minimum height=0.5cm, above left] {$m+1$};
\draw (-1.25,0.5) node [minimum width=0.5cm, minimum height=0.5cm, above left] {$m$};
\end{tikzpicture}
\end{figure}
so that we can apply  chute move-1 (here numbers, located under the boxes of the $i$-th row
and below the boxes of the $i+1$-th row are the indexes of the corresponding antidiagonals and the letters $\emph{A, B,
A}'$ and $\emph{B}'$ are the corresponding subwords, obtained by "reading"  of the pipe dream $D$ in the way, described above).
Then the subword obtained by  reading the $m$-th and $(m+1)$-th rows is the following:

$$As_{i+k}s_{i+k-1}\ldots s_{i+1}s_{i}A'Bs_{i+k}s_{i+k-1}\ldots s_{i+2}s_{i+1}B'.$$

From the properties of the Demazure product it follows that this subword is equivalent to the following:

$$AA's_{i+k}s_{i+k-1}\ldots s_{i+1}s_{i}s_{i+k}s_{i+k-1}\ldots s_{i+2}s_{i+1}BB'.$$

This means that for the case of chute move-1 we are reduced to proving the equivalence of the words
$$(s_{i+k}s_{i+k-1}\ldots s_{i+1}s_{i}s_{i+k}s_{i+k-1}\ldots s_{i+2}s_{i+1})\ \textmd{and}\
(s_{i+k-1}s_{i+k-2}\ldots s_{i+1}s_{i}s_{i+k}s_{i+k-1}\ldots s_{i+1}s_{i})$$
with respect to the Demazure product. We will carry out the proof inductively (the parameter
of induction is  the $\emph{move length}$ denoted by $L$, by which we mean the length of the
corresponding $2\times L$ $\emph{chutable rectangle}$).
The inductive basis is the simple equivalence of the words $(s_{i}s_{i+1}s_{i})\cong(s_{i+1}s_{i}s_{i+1})$.

Let us suppose that our statement is proven for the case $L\leq k+1$. Then for our
initial word $(s_{i+k}s_{i+k-1}\ldots s_{i+1}s_{i}s_{i+k}s_{i+k-1}\ldots s_{i+2}s_{i+1)}$ the following
equivalence relations take place:

$$(s_{i+k}s_{i+k-1}\ldots s_{i+1}s_{i}\widehat{s_{i+k}}s_{i+k-1}\ldots s_{i+2}s_{i+1})\cong$$

$$(\textrm{obtained by shifting} \ s_{i+k} \ \textrm{to the left until "collision"  with } \ s_{i+k})$$

$$\cong(s_{i+k}s_{i+k-1}s_{i+k})(s_{i+k-2}\ldots s_{i}s_{i+k-1}\ldots s_{i+1})\cong$$

$$(\textrm{obtained by applying the equivalence relation} \ s_{i+k}s_{i+k-1}s_{i+k}\rightarrow s_{i+k-1}s_{i+k}s_{i+k-1})$$

$$\cong(s_{i+k-1}s_{i+k})(s_{i+k-1}\ldots s_{i}s_{i+k-1}\ldots s_{i+1})\cong$$

$$(\textrm{obtained by appying the induction hypothesis})$$

$$\cong(s_{i+k-1}\widehat{s_{i+k}})(s_{i+k-2}\ldots s_{i}s_{i+k-1}\ldots s_{i})\cong$$

$$(\textrm{obtained by shifting} \ s_{i+k} \ \textrm{to the right until "collision"  with} \ s_{i+k-1})$$

$$\cong(s_{i+k-1}s_{i+k-2}\ldots s_{i+1}s_{i}s_{i+k}\ldots s_{i+1}s_{i}).$$

So the statement is proven for the case $L=k+2$ and the induction step is fulfilled.

2) Let us now suppose that the $m$-th and $(m+1)$-th rows of $D$ look like the figure above,

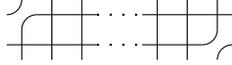
\begin{figure}[!ht]
\centering
\begin{tikzpicture}
\draw[step=0.5cm]  (0,0) grid (1.5,1);
\draw[xshift=2.25cm, step=0.5cm]  (0,0) grid (1.5,1);
\draw (1.5,1)--(2.25,1);
\draw (1.5,0)--(2.25,0);
\draw[dotted, thick] (1.5,0.5)--(2.25,0.5);
\draw (0,0) node [minimum width=0.5cm, minimum height=0.5cm, above left] {$A'$};
\draw (0,0.5) node [minimum width=0.5cm, minimum height=0.5cm, above left] {$B'$};
\draw (3.75,0) node [minimum width=0.5cm, minimum height=0.5cm, above right] {$B$};
\draw (3.75,0.5) node [minimum width=0.5cm, minimum height=0.5cm, above right] {$A$};

\draw (0.5,0) node [minimum width=0.5cm, minimum height=0.5cm, below left] {$\scriptscriptstyle(i)$};
\draw (0.33,0) node [minimum width=0.5cm, minimum height=0.5cm, below right] {$\scriptscriptstyle(i+1)$};

\draw (1,1) node [minimum width=0.5cm, minimum height=0.5cm, above left] {$\scriptscriptstyle(i)$};
\draw (0.83,1) node [minimum width=0.5cm, minimum height=0.5cm, above right] {$\scriptscriptstyle(i+1)$};

\draw (2.58,0) node [minimum width=0.5cm, minimum height=0.5cm, below right] {$\scriptscriptstyle(i+k)$};
\draw (3.08,1) node [minimum width=0.5cm, minimum height=0.5cm, above right] {$\scriptscriptstyle(i+k)$};

\draw (0.5,0.5) node [minimum width=0.5cm, minimum height=0.5cm, above right] {$+$};
\draw (0.5,0) node [minimum width=0.5cm, minimum height=0.5cm, above right] {$+$};
\draw (1,0.5) node [minimum width=0.5cm, minimum height=0.5cm, above right] {$+$};
\draw (1,0) node [minimum width=0.5cm, minimum height=0.5cm, above right] {$+$};

\draw[xshift=1.75cm] (0.5,0.5) node [minimum width=0.5cm, minimum height=0.5cm, above right] {$+$};
\draw[xshift=1.75cm] (0.5,0) node [minimum width=0.5cm, minimum height=0.5cm, above right] {$+$};
\draw[xshift=1.75cm] (1,0.5) node [minimum width=0.5cm, minimum height=0.5cm, above right] {$+$};
\draw[xshift=1.75cm] (1,0) node [minimum width=0.5cm, minimum height=0.5cm, above right] {$+$};
\draw[xshift=2.75cm] (0.5,0.5) node [minimum width=0.5cm, minimum height=0.5cm, above right] {$+$};
\draw (0,0) node [minimum width=0.5cm, minimum height=0.5cm, above right] {$+$};

\draw (-1,0) node [minimum width=0.5cm, minimum height=0.5cm, above left] {$m+1$};
\draw (-1.25,0.5) node [minimum width=0.5cm, minimum height=0.5cm, above left] {$m$};
\end{tikzpicture}
\end{figure}

so that we can apply chute move-2. By carrying out the arguments analogously to 1), we
reduce our statement  to proving the equivalence of the words
$$(s_{i+k}s_{i+k-1}\ldots s_{i+1}s_{i}s_{i+k}s_{i+k-1}\ldots s_{i+1}s_{i})\ \textmd{and}\
(s_{i+k-1}s_{i+k-2}\ldots s_{i+1}s_{i}s_{i+k}s_{i+k-1}\ldots s_{i+1}s_{i}).$$

Again we apply the induction on the parameter $L$. The induction basis is now the equivalence relation
$(s_{i}s_{i})\cong s_{i}$. Then, analogously to 1), for our initial word
$(s_{i+k}s_{i+k-1}\ldots s_{i+1}$ $s_{i}s_{i+k}s_{i+k-1}\ldots s_{i+1}s_{i})$ the following equivalence
relations take place:

$$(s_{i+k}s_{i+k-1}\ldots s_{i+1}\widehat{s_{i}}s_{i+k}s_{i+k-1}\ldots s_{i+1}s_{i})\cong$$

$$(\textrm{obtained by shifting} \ s_{i}  \ \textrm{to the right until "collision" with} \ s_{i+1})$$

$$\cong(s_{i+k}s_{i+k-1}\ldots s_{i+1}s_{i+k}s_{i+k-1}\ldots s_{i+2})(s_{i}s_{i+1}s_{i})\cong$$

$$(\textrm{obtained by applying the equivalence relation} \ s_{i}s_{i+1}s_{i}\rightarrow s_{i+1}s_{i}s_{i+1})$$

$$\cong(s_{i+k}s_{i+k-1}\ldots s_{i+1}s_{i+k}s_{i+k-1}\ldots s_{i+1})(s_{i}s_{i+1})\cong$$

$$(\textrm{obtained by appying the induction hypothesis})$$

$$\cong(s_{i+k-1}s_{i+k-2}\ldots s_{i+1}s_{i+k}s_{i+k-1}\ldots s_{i+1})(s_{i}s_{i+1})\cong$$

$$(\textrm{obtained by applying the equivalence relation} \ s_{i+1}s_{i}s_{i+1}\rightarrow s_{i}s_{i+1}s_{i})$$

$$\cong(s_{i+k-1}s_{i+k-2}\ldots s_{i+1}s_{i+k}s_{i+k-1}\ldots\widehat{s_{i}})(s_{i+1}s_{i})\cong$$

$$(\textrm{obtained by shifting} \ s_{i}  \ \textrm{to the left until "collision" until} \ s_{i+1})$$

$$\cong(s_{i+k-1}s_{i+k-2}\ldots s_{i+1}s_{i}s_{i+k}s_{i+k-1}\ldots s_{i+1}s_{i}),$$

Q.E.D.$\blacksquare$

$\textbf{Remark 2.1:}$ Note that in case of reduced pipe dreams chute moves  2 and 3 cannot be used,
because their chutable rectangles are given by  double crossing of the pair of strands:
\begin{figure}[!ht]
\centering
\begin{tikzpicture}[scale=0.4]
\draw (0,0.5)--(3,0.5);
\draw (1,1.5)--(3,1.5);
\draw (1.5,0)--(1.5,2);
\draw (2.5,0)--(2.5,2);
\draw (0.5,0)--(0.5,1);

\draw[loosely dotted, thick] (3,0.5)--(4.5,0.5);
\draw[loosely dotted, thick] (3,1.5)--(4.5,1.5);

\draw (4.5,0.5)--(6.5,0.5);
\draw (4.5,1.5)--(7.5,1.5);
\draw (5,0)--(5,2);
\draw (6,0)--(6,2);
\draw (7,1)--(7,2);

\draw (1,1.5) arc (90:180:0.5cm);
\draw (0,1.5) arc (270:360:0.5cm);
\draw (7.5,0.5) arc (90:180:0.5cm);
\draw (6.5,0.5) arc (270:360:0.5cm);
\end{tikzpicture}
\end{figure}

$$$$
$$$$
Also in case of applying chute move-1 to a reduced pipe dream preservation of the permutation can easily be proven
graphically. Indeed,  all stands, except the pair involved in the conversion, remain "untouched".
Also the crossing of two involved strands in the upper right corner of the chutable rectangle is replaced with
the crossing in the lower left corner and the exits of these two strands stay the same (see
the figure below).
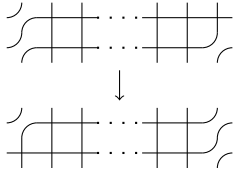
\begin{figure}[!ht]
\centering
\begin{tikzpicture}[scale=0.4]
\draw (1,0.5)--(3,0.5);
\draw (1,1.5)--(3,1.5);
\draw (1.5,0)--(1.5,2);
\draw (2.5,0)--(2.5,2);
\draw (1,0.5) arc (90:180:0.5cm);
\draw (0,0.5) arc (270:360:0.5cm);

\draw[loosely dotted, thick] (3,0.5)--(4.5,0.5);
\draw[loosely dotted, thick] (3,1.5)--(4.5,1.5);

\draw (4.5,0.5)--(6.5,0.5);
\draw (4.5,1.5)--(7.5,1.5);
\draw (5,0)--(5,2);
\draw (6,0)--(6,2);
\draw (7,1)--(7,2);

\draw (1,1.5) arc (90:180:0.5cm);
\draw (0,1.5) arc (270:360:0.5cm);
\draw (7.5,0.5) arc (90:180:0.5cm);
\draw (6.5,0.5) arc (270:360:0.5cm);

\draw[yshift=-3.5cm] (1,1.5)--(3,1.5);
\draw[yshift=-3.5cm] (0,0.5)--(3,0.5);
\draw[yshift=-3.5cm] (1.5,0)--(1.5,2);
\draw[yshift=-3.5cm] (2.5,0)--(2.5,2);
\draw[yshift=-3.5cm] (7.5,1.5) arc (90:180:0.5cm);
\draw[yshift=-3.5cm] (6.5,1.5) arc (270:360:0.5cm);

\draw[yshift=-3.5cm, loosely dotted, thick] (3,0.5)--(4.5,0.5);
\draw[yshift=-3.5cm,loosely dotted, thick] (3,1.5)--(4.5,1.5);

\draw[yshift=-3.5cm] (4.5,0.5)--(6.5,0.5);
\draw[yshift=-3.5cm] (4.5,1.5)--(6.5,1.5);
\draw[yshift=-3.5cm] (5,0)--(5,2);
\draw[yshift=-3.5cm] (6,0)--(6,2);
\draw[yshift=-3.5cm] (0.5,0)--(0.5,1);

\draw[yshift=-3.5cm] (1,1.5) arc (90:180:0.5cm);
\draw[yshift=-3.5cm] (0,1.5) arc (270:360:0.5cm);
\draw[yshift=-3.5cm] (7.5,0.5) arc (90:180:0.5cm);
\draw[yshift=-3.5cm] (6.5,0.5) arc (270:360:0.5cm);

\draw [->] (3.75,-0.25)--(3.75,-1.25);
\end{tikzpicture}
\end{figure}

Now we can redefine the mitosis algorithm  in terms of chute moves:

$\textbf{Proposition 2.3:}$ Let $D$ be a pipe dream and $j_{min}$ be a minimal column index, such that
$(i+1,j)\notin D$ and for any $p\leqslant j$ $(i,p)$ belongs to $D$. Then each $D_{p}\in mitosis_{i}(D)$
obtained from $D$ by

1) deleting $(i,j_{min})$, and then,

2) moving to the right, one by one applying chute moves-1, so that  $(i,p)$ is the last cross,
moved from the $i$-th row to $i+1$-th row.

From this proposition we can see that all elements of $mitosis_{i}(D)$ correspond to the same permutation,
because all of them are the results of chute moves, applied to $D\backslash(i,j_{min})$. Now we are ready
to prove the main statement of this section:

$\textbf{Theorem 2.3:}$ Let $w$ be an element of $S_{\infty}$, $s_{i}$ --- an elementary transposition and let 
$D$ belong to $\mathcal{P}(w)$,
with the condition that $l(ws_{i})=l(w)-1$. Then, if the set $mitosis_{i}(D)$ is not empty, then it
lies entirely in either $\mathcal{P}(ws_{i})$ or $\mathcal{P}(w)$.

$\textbf{Proof:}$ From our reasoning above we can see that if $D\backslash(i,j_{min})$ belongs to $\mathcal{P}(ws_{i})$
(accordingly, to $\mathcal{P}(w)$),then the same is true of all the elements of $mitosis_{i}(D)$.

Let $D$ look like
$$$$
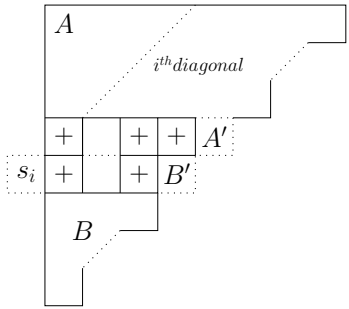
\begin{figure}
\centering
\scalebox{0.5}{
\begin{tikzpicture}
\draw (1,1)--(1,0)--(0,0)--(0,8)--(8,8)--(8,7)--(7,7);
\draw[loosely dotted, thick] (7,7)--(6,6);
\draw (6,6)--(6,5)--(5,5);
\draw[loosely dotted, thick] (5,5)--(5,4)--(4,4)--(4,3)--(3,3);
\draw (3,3)--(3,2)--(2,2);
\draw[loosely dotted, thick] (2,2)--(1,1);
\draw[loosely dotted, thick] (0,4)--(-1,4)--(-1,3)--(0,3);
\draw[font=\huge] (-1,3) node [minimum width=1cm, minimum height=1cm, above right] {$s_{i}$};
\draw (0,3)--(3,3)--(3,5)--(0,5);
\draw (2,4)--(4,4)--(4,5)--(3,5);
\draw (0,4)--(1,4)--(1,5)--(1,3);
\draw[loosely dotted, thick] (1,4)--(2,4);
\draw[loosely dotted, thick] (4,5)--(5,5);
\draw (2,3)--(2,5);

\draw[font=\huge] (0,3) node [minimum width=1cm, minimum height=1cm, above right] {$+$};
\draw[font=\huge] (0,4) node [minimum width=1cm, minimum height=1cm, above right] {$+$};
\draw[font=\huge] (2,3) node [minimum width=1cm, minimum height=1cm, above right] {$+$};
\draw[font=\huge] (2,4) node [minimum width=1cm, minimum height=1cm, above right] {$+$};
\draw[font=\huge] (3,4) node [minimum width=1cm, minimum height=1cm, above right] {$+$};
\draw[font=\huge] (4,4) node [minimum width=1cm, minimum height=1cm, above right] {$A'$};
\draw[font=\huge] (3,3) node [minimum width=1cm, minimum height=1cm, above right] {$B'$};
\draw[font=\huge] (0.5,1.5) node [minimum width=1cm, minimum height=1cm, above right] {$B$};
\draw[font=\huge] (0,7) node [minimum width=1cm, minimum height=1cm, above right] {$A$};
\draw[loosely dotted,thick] (1,5)--(4,8);
\draw[font=\Large] (2.75,6) node [above right] {$i^{th} diagonal$};
\end{tikzpicture}
}
\end{figure}
$$$$

For every pipe dream $B$  the corresponding word will be denoted by $word(B)$. For every word
$v$ and pipe dream $B$ the corresponding Demazure product will be denoted by $Demaz(v)$ ($Demaz(B)$
accordingly). Then the following equalities hold:

$$w=Demaz(D)=Demaz(AA's_{i+j_{min}-1}\ldots s_{i}B's_{i+j_{min}-1}\ldots s_{i+1}B).$$

By applying to the left word a transform, analogous to chute move-1, and given that $s_{i}$ commutes with
the subword $B$, obtain

$$w=Demaz(D)=Demaz(AA's_{i+j_{min}-2}\ldots s_{i}B's_{i+j_{min}-1}\ldots s_{i+1}Bs_{i})=$$
$$=Demaz(word(D')s_{i}).$$

Here $D'$ is the result of removing $(i,j_{min})$ from $D$.

Denote by $\tilde{w}$ the permutation $Demaz(word(D'))$. Then, according to the definition
of the Demazure product, two cases are possible:

1)$l(\tilde{w}s_{i})=l(\tilde{w})-1$. Then we have $w=Demaz(word(D')s_{i})=\tilde{w}$
and $D'$ belongs to $\mathcal{P}(w)$.

2)$l(\tilde{w}s_{i})=l(\tilde{w})+1$. Then we have $w=Demaz(word(D')s_{i}))=\tilde{w}s_{i}$,
which means that $\tilde{w}=ws_{i}$ and $D'$ belongs to $\mathcal{P}(ws_{i})$ (note, that in both
cases $l(ws_{i})=l(w)-1$). Q.E.D. $\blacksquare$

In compliance with the proven theorem we divide $\mathcal{P}(w)$ into three
disjoint sets: $\mathcal{P}_{s}(w)$ (the set of all pipe dreams,which are sent to $\mathcal{P}(ws_{i})$),
$\mathcal{P}_{I}(w)$ (the set of all pipe dreams, which are sent to $\mathcal{P}(w)$), and
$\mathcal{P}_{\varnothing}(w)$ (the set of all pipe dreams, which are sent to to the empty set). Here and further
by mitosis we mean $mitosis_{i}$ with condition that $l(ws_{i})=l(w)-1$. This partition
will be used later.

$$$$
%%%%%%%%%%%%%%%%%%%%%%%%%%%%%%%%%%%%%%%%%%%%%%%%%%%%%%%%%%%%%%%%%%%%%%%%%
\section{Intron mutations}%%%%%%%%%%%%%%%%%%%%%%%%%%%%%%%%%%%%%%%%%%%%%%%%%%%%
%%%%%%%%%%%%%%%%%%%%%%%%%%%%%%%%%%%%%%%%%%%%%%%%%%%%%%%%%%%%%%%%%%%%%%%%%
$$$$

Let $D$ be an arbitrary pipe dream with a fixed row index $i$. Index the boxes in
in $i$-th and $(i+1)$-th rows as shown in the following figure:
\begin{figure}[!ht]
\centering
\scalebox{0.5}{
\begin{tikzpicture}
\draw (0,0) grid (4,2);
\draw (8,0)--(4,0)--(4,2)--(8,2);
\draw[loosely dotted, very thick] (4.5,1)--(6,1);
\draw[font=\huge] (-0.5,0) node [minimum width=1cm, minimum height=1cm, above left] {$(i+1)$};
\draw[font=\huge] (-1.15,1) node [minimum width=1cm, minimum height=1cm, above left] {$(i)$};
\draw[font=\huge] (0,1) node [minimum width=1cm, minimum height=1cm, above right] {$1$}; 
\draw[font=\huge] (0,0) node [minimum width=1cm, minimum height=1cm, above right] {$2$};
\draw[font=\huge] (1,1) node [minimum width=1cm, minimum height=1cm, above right] {$3$}; 
\draw[font=\huge] (1,0) node [minimum width=1cm, minimum height=1cm, above right] {$4$};
\draw[font=\huge] (2,1) node [minimum width=1cm, minimum height=1cm, above right] {$5$}; 
\draw[font=\huge] (2,0) node [minimum width=1cm, minimum height=1cm, above right] {$6$}; 
\draw[font=\huge] (3,1) node [minimum width=1cm, minimum height=1cm, above right] {$7$};  
\draw[font=\huge] (3,0) node [minimum width=1cm, minimum height=1cm, above right] {$8$}; 
\draw[font=\huge] (0,2) node [minimum width=1cm, minimum height=1cm, above right] {$1$};  
\draw[font=\huge] (1,2) node [minimum width=1cm, minimum height=1cm, above right] {$2$}; 
\draw[font=\huge] (2,2) node [minimum width=1cm, minimum height=1cm, above right] {$3$};  
\draw[font=\huge] (3,2) node [minimum width=1cm, minimum height=1cm, above right] {$4$}; 
\draw[loosely dotted, very thick] (4.5,2.5)--(6,2.5); 
\end{tikzpicture}
}
\end{figure}

Then henceforth by an $\emph{intron}$ we mean a $2\times m$ rectangle, located in
these two adjasent rows such that:

1) the first and the last boxes of this rectangle are empty,

2) no $\elpl$ column can be located to the right from a $\elel$ column
or a $\plel$ column and no $\elel$ column can be located to the right from a $\plel$ column.
(The first (last) box of a rectangle is the box with the maximal (minimal) index
according to the ordering described above.)

An intron $C$ is $\emph{maximal}$, if the empty box with largest index before $C$ (if there is one) 
is located in the $i+1$-th row and the empty box with smallest index after $C$,
is located in the $i$-th row. In other words, an intron is maximal, if in cannot
be extended rightwards or leftwards.

$\textbf{Lemma 3.1:}$ Let $D$ be a pipe dream and $C\subseteq D$ be an intron. Then by applying a sequence of chute moves
and inverse chute moves we can transform $C$ to a new intron $\tau(C)$ with the following properties:

1) the set of $\plpl$ columns in $C$ coincide with the set of $\plpl$ columns in $\tau(C)$, and

2) the number $c_{i}$ of crosses in the $i$-row of $C$ coincide with the number $\tilde{c}_{i+1}$
of crosses in the $i+1$-th row of $\tau(C)$ and vise versa.

$\textbf{Proof:}$ Suppose that $c_{i}\geq\tilde{c}_{i+1}$. Then we will carry out the proof
inductively with the parameter $c=c_{i}-\tilde{c}_{i+1}$. In case $c=0$ we obviously have
$C=\tau(C)$. Now if $c>0$ we consider the leftmost $\plel$ column (denote its index
by $p$). Moving to the left from this column we will, sooner of later, find a column 
of the type $\elel$ or $\elpl$. If it is a $\elel$ column, then we can apply chute move-1 and thereby
chute the cross from the $i$-th row to the $i+1$-th one. Since the result of this conversion will also 
be an intron, the proof is reduced to the induction hypothesis. If it looks like a $\elpl$, then, owing to 
the fact that $c>0$ and, thereafter, there is more than one $\plel$ column in $C$, the corresponding
fragment of $C$ will look like

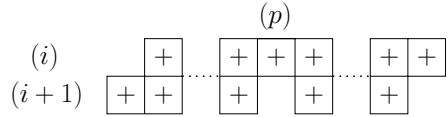
\begin{figure}[!ht]
\centering
\scalebox{0.5}{
\begin{tikzpicture}
\draw (2,1)--(0,1)--(0,0)--(2,0)--(2,2)--(1,2)--(1,0);
\draw (4,2)--(4,0)--(3,0)--(3,2)--(6,2)--(6,0)--(5,0)--(5,2)--(4,2);
\draw (3,1)--(6,1);
\draw (8,2)--(8,0)--(7,0)--(7,2)--(9,2)--(9,1)--(7,1);
\draw[loosely dotted, very thick] (2,1)--(3,1);
\draw[loosely dotted, very thick] (6,1)--(7,1);
\draw[font=\huge] (-0.5,0) node [minimum width=1cm, minimum height=1cm, above left] {$(i+1)$};
\draw[font=\huge] (-1.15,1) node [minimum width=1cm, minimum height=1cm, above left] {$(i)$}; 
\draw[font=\huge] (3.9,2.1) node [minimum width=1cm, minimum height=1cm, above right] {$(p)$}; 
\draw[font=\huge] (0,0) node [minimum width=1cm, minimum height=1cm, above right] {$+$};  
\draw[font=\huge] (1,1) node [minimum width=1cm, minimum height=1cm, above right] {$+$}; 
\draw[font=\huge] (1,0) node [minimum width=1cm, minimum height=1cm, above right] {$+$}; 
\draw[font=\huge] (3,0) node [minimum width=1cm, minimum height=1cm, above right] {$+$};
\draw[font=\huge] (3,1) node [minimum width=1cm, minimum height=1cm, above right] {$+$};
\draw[font=\huge] (4,1) node [minimum width=1cm, minimum height=1cm, above right] {$+$};
\draw[font=\huge] (5,0) node [minimum width=1cm, minimum height=1cm, above right] {$+$};
\draw[font=\huge] (5,1) node [minimum width=1cm, minimum height=1cm, above right] {$+$};
\draw[font=\huge] (7,0) node [minimum width=1cm, minimum height=1cm, above right] {$+$};
\draw[font=\huge] (7,1) node [minimum width=1cm, minimum height=1cm, above right] {$+$};
\draw[font=\huge] (8,1) node [minimum width=1cm, minimum height=1cm, above right] {$+$};
\end{tikzpicture}

}
\end{figure}

Then by applying  chute move-2 and reverse chute move-3 the way it shown on the following figure, we will
bring the corresponding fragment of $C$ to a form

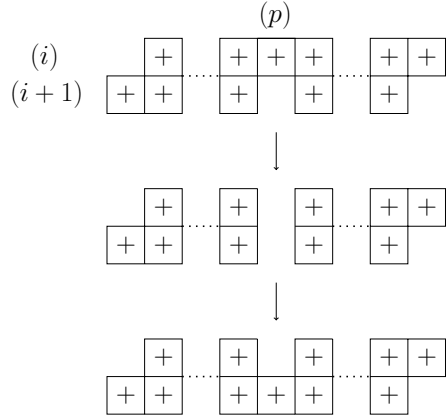
\begin{figure}[!ht]
\centering
\scalebox{0.5}{
\begin{tikzpicture}
\draw (2,1)--(0,1)--(0,0)--(2,0)--(2,2)--(1,2)--(1,0);
\draw (4,2)--(4,0)--(3,0)--(3,2)--(6,2)--(6,0)--(5,0)--(5,2)--(4,2);
\draw (3,1)--(6,1);
\draw (8,2)--(8,0)--(7,0)--(7,2)--(9,2)--(9,1)--(7,1);
\draw[loosely dotted, very thick] (2,1)--(3,1);
\draw[loosely dotted, very thick] (6,1)--(7,1);
\draw[font=\huge] (-0.5,0) node [minimum width=1cm, minimum height=1cm, above left] {$(i+1)$};
\draw[font=\huge] (-1.15,1) node [minimum width=1cm, minimum height=1cm, above left] {$(i)$}; 
\draw[font=\huge] (3.9,2.1) node [minimum width=1cm, minimum height=1cm, above right] {$(p)$}; 
\draw[font=\huge] (0,0) node [minimum width=1cm, minimum height=1cm, above right] {$+$};  
\draw[font=\huge] (1,1) node [minimum width=1cm, minimum height=1cm, above right] {$+$}; 
\draw[font=\huge] (1,0) node [minimum width=1cm, minimum height=1cm, above right] {$+$}; 
\draw[font=\huge] (3,0) node [minimum width=1cm, minimum height=1cm, above right] {$+$};
\draw[font=\huge] (3,1) node [minimum width=1cm, minimum height=1cm, above right] {$+$};
\draw[font=\huge] (4,1) node [minimum width=1cm, minimum height=1cm, above right] {$+$};
\draw[font=\huge] (5,0) node [minimum width=1cm, minimum height=1cm, above right] {$+$};
\draw[font=\huge] (5,1) node [minimum width=1cm, minimum height=1cm, above right] {$+$};
\draw[font=\huge] (7,0) node [minimum width=1cm, minimum height=1cm, above right] {$+$};
\draw[font=\huge] (7,1) node [minimum width=1cm, minimum height=1cm, above right] {$+$};
\draw[font=\huge] (8,1) node [minimum width=1cm, minimum height=1cm, above right] {$+$};
\begin{scope}[yshift=-4cm]
\draw (2,1)--(0,1)--(0,0)--(2,0)--(2,2)--(1,2)--(1,0);
\draw (4,2)--(4,0)--(3,0)--(3,2)--(4,2);
\draw (6,2)--(6,0)--(5,0)--(5,2)--(6,2);
\draw (3,1)--(4,1);
\draw (5,1)--(6,1);
\draw (8,2)--(8,0)--(7,0)--(7,2)--(9,2)--(9,1)--(7,1);
\draw[loosely dotted, very thick] (2,1)--(3,1);
\draw[loosely dotted, very thick] (6,1)--(7,1);  
\draw[font=\huge] (0,0) node [minimum width=1cm, minimum height=1cm, above right] {$+$};  
\draw[font=\huge] (1,1) node [minimum width=1cm, minimum height=1cm, above right] {$+$}; 
\draw[font=\huge] (1,0) node [minimum width=1cm, minimum height=1cm, above right] {$+$}; 
\draw[font=\huge] (3,0) node [minimum width=1cm, minimum height=1cm, above right] {$+$};
\draw[font=\huge] (3,1) node [minimum width=1cm, minimum height=1cm, above right] {$+$};
\draw[font=\huge] (5,0) node [minimum width=1cm, minimum height=1cm, above right] {$+$};
\draw[font=\huge] (5,1) node [minimum width=1cm, minimum height=1cm, above right] {$+$};
\draw[font=\huge] (7,0) node [minimum width=1cm, minimum height=1cm, above right] {$+$};
\draw[font=\huge] (7,1) node [minimum width=1cm, minimum height=1cm, above right] {$+$};
\draw[font=\huge] (8,1) node [minimum width=1cm, minimum height=1cm, above right] {$+$};
\end{scope}

\begin{scope}[yshift=-8cm]
\draw (2,1)--(0,1)--(0,0)--(2,0)--(2,2)--(1,2)--(1,0);
\draw (4,2)--(4,0)--(3,0)--(3,2)--(4,2);
\draw (6,2)--(6,0)--(5,0)--(5,2)--(6,2);
\draw (3,1)--(6,1);
\draw (4,0)--(5,0);
\draw (8,2)--(8,0)--(7,0)--(7,2)--(9,2)--(9,1)--(7,1);
\draw[loosely dotted, very thick] (2,1)--(3,1);
\draw[loosely dotted, very thick] (6,1)--(7,1);  
\draw[font=\huge] (0,0) node [minimum width=1cm, minimum height=1cm, above right] {$+$};  
\draw[font=\huge] (1,1) node [minimum width=1cm, minimum height=1cm, above right] {$+$}; 
\draw[font=\huge] (1,0) node [minimum width=1cm, minimum height=1cm, above right] {$+$}; 
\draw[font=\huge] (3,0) node [minimum width=1cm, minimum height=1cm, above right] {$+$};
\draw[font=\huge] (3,1) node [minimum width=1cm, minimum height=1cm, above right] {$+$};
\draw[font=\huge] (5,0) node [minimum width=1cm, minimum height=1cm, above right] {$+$};
\draw[font=\huge] (5,1) node [minimum width=1cm, minimum height=1cm, above right] {$+$};
\draw[font=\huge] (7,0) node [minimum width=1cm, minimum height=1cm, above right] {$+$};
\draw[font=\huge] (7,1) node [minimum width=1cm, minimum height=1cm, above right] {$+$};
\draw[font=\huge] (8,1) node [minimum width=1cm, minimum height=1cm, above right] {$+$};
\draw[font=\huge] (4,0) node [minimum width=1cm, minimum height=1cm, above right] {$+$};
\end{scope}
\draw[->] (4.5,-0.5)--(4.5,-1.5);
\draw[yshift=-4cm, ->] (4.5,-0.5)--(4.5,-1.5);

\end{tikzpicture}

}
\end{figure}
and thereby again chute the cross from the $i$-th row to the $i+1$-th one. Thus the proof
is again reduced to the induction hypothesis.

In case $c_{i}<\tilde{c}_{i+1}$ we just flip the argument $180\textdegree$. $\blacksquare$

The transformation $\tau$ is called $\emph{intron mutation}$. Note, that the intron $\tau(C)$ is defined uniquely
and, by constraction, $\tau(\tau(C))=C$, i.e. $\tau$ is an involution.

Now we are ready to the main statement of this section:

$\textbf{Theorem 3.1:}$ Let $w$ be an element of $S_{\infty}$. Then for each $i\in\mathbb{N}$
there is an involution $\tau_{i}:\mathcal{P}(w)\longrightarrow\mathcal{P}(w)$, such that for any
$D\in\mathcal{P}(w)$ the following conditions take place:

1) $\tau_{i}(D)$ coincides with $D$ in all rows with indexes, different from $i$ and $i+1$.

2) $start_{i}(D)=start_{i}(\tau_{i}(D))$ and $\tau_{i}(D)$ agrees with $D$ in all columns
with indices strictly less than $start_{i}(D)$.

3) $l_{i}^{i}(\tau_{i}(D))=l_{i+1}^{i}(D)$ (here $l_{r}^{i}(-)$ - is the number of crosses in the $r$-th row,
located to the right or in column with index $start_{i}(-)$)

$\textbf{Proof:}$ Let $D$ belong to $\mathcal{P}(w)$. Consider all crosses of the union of the $i$-th 
and $(i+1)$-th rows, located to the right or in column with index $start_{i}(D)$. Then, according to 
the definition of  $start_{i}(D)$, we can find a minimal rectangle with empty last box,
starting at the column $start_{i}(D)$ and containing all these crosses. 
Since the first box of this rectangle is also empty, it can be uniquely represented in form of
a disjoint union of maximal introns and rectangles, completely filled with crosses. Apply to an
every maximal intron the intron mutation, described above. Since every mutation is a sequence
of chute moves and reverse chute moves, the obtained pipe dream obviously belongs to $\mathcal{P}(w)$. 
Owing the fact that each mutation is an involution and that the result of applying a mutation to a maximal
intron is also a maximal intron, the obtained transformation is also an involution. Properties 1)-3) are obvious
from the consctruction scheme. $\blacksquare$

$\textbf{Remark 3.1:}$ Note that the partition 
$\mathcal{P}(w)=\mathcal{P}_{\varnothing}(w)\sqcup\mathcal{P}_{s}(w)\sqcup\mathcal{P}_{I}(w)$ 
is invariant under the constructed involution. This fact will be used henceforth.

$$$$
%%%%%%%%%%%%%%%%%%%%%%%%%%%%%%%%%%%%%%%%%%%%%%%%%%%%%%%%%%%%%%%%%%%%%%%%%
\section{Mitosis theorem}%%%%%%%%%%%%%%%%%%%%%%%%%%%%%%%%%%%%%%%%%%%%%%%%%%%%
%%%%%%%%%%%%%%%%%%%%%%%%%%%%%%%%%%%%%%%%%%%%%%%%%%%%%%%%%%%%%%%%%%%%%%%%%
$$$$

In the second section we have defined the way mitosis acts on the set $\mathcal{P}(w)$. 
Nevertheless, in order to prove the main theorem of this article we will have to slightly
modify it's initial definition:

$\textbf{Definition 4.1:}$ Let $l(ws_{i})=l(w)-1$, $D$ belong to $\mathcal{P}(w)$ and
$\mathcal{J}_{i}(D)={j_{1},\ldots,j_{k}}$. Then operator $mitosis'_{i}$ sends $D$ to the set

$$\{D_{j_{1}}; D_{j_{1}}+D_{j_{2}}; D_{j_{2}};\ldots;D_{j_{k-1}};D_{j_{k-1}}+D_{j_{k}};D_{j_{k}}\}.$$

Here by "sum" we mean the union of the elements of the corresponding pipe dreams.

(As we can see, the elements of the set $mitosis'_{i}(D)$ form some kind of a chain, where links
are elements of $D_{j_{m}}$ and the result of the two links’ cohesion is a sum of the corresponding pipe dreams).

Let us show that pipe dreams of $mitosis'_{i}(D)$ represent the same permutation.
Indeed, let $D_{j_{m}}$ and $D_{j_{m+1}}$ look like
\begin{figure}[!ht]
\centering
\scalebox{0.6}{
\begin{tikzpicture}
 
\draw (0,0) grid (2,2);
\draw (3,0) grid (5,2);

\draw[font=\huge] (0,0.5) node [minimum width=1cm, minimum height=1cm, above left] {$\cdots$};
\draw[font=\huge] (5,0.5) node [minimum width=1cm, minimum height=1cm, above right] {$\cdots$}; 
\draw[font=\huge] (2,0.5) node [minimum width=1cm, minimum height=1cm, above right] {$\cdots$}; 
\draw[font=\huge] (-1,0.5) node [minimum width=1cm, minimum height=1cm, above left] {$D_{j_{m}}$}; 
\draw[font=\huge] (1,0) node [minimum width=1cm, minimum height=1cm, above right] {$+$}; 
\draw[font=\huge] (1,1) node [minimum width=1cm, minimum height=1cm, above right] {$+$}; 
\draw[font=\huge] (3,0) node [minimum width=1cm, minimum height=1cm, above right] {$+$}; 
\draw[font=\huge] (3,1) node [minimum width=1cm, minimum height=1cm, above right] {$+$}; 
\draw[font=\huge] (4,1) node [minimum width=1cm, minimum height=1cm, above right] {$+$};

\begin{scope}[yshift=-3cm]
\draw (0,0) grid (2,2);
\draw (3,0) grid (5,2);

\draw[font=\huge] (0,0.5) node [minimum width=1cm, minimum height=1cm, above left] {$\cdots$};
\draw[font=\huge] (5,0.5) node [minimum width=1cm, minimum height=1cm, above right] {$\cdots$}; 
\draw[font=\huge] (2,0.5) node [minimum width=1cm, minimum height=1cm, above right] {$\cdots$}; 
\draw[font=\huge] (-1,0.5) node [minimum width=1cm, minimum height=1cm, above left] {$D_{j_{m+1}}$}; 
\draw[font=\huge] (1,0) node [minimum width=1cm, minimum height=1cm, above right] {$+$}; 
\draw[font=\huge] (1,1) node [minimum width=1cm, minimum height=1cm, above right] {$+$}; 
\draw[font=\huge] (3,0) node [minimum width=1cm, minimum height=1cm, above right] {$+$}; 
\draw[font=\huge] (3,1) node [minimum width=1cm, minimum height=1cm, above right] {$+$}; 
\draw[font=\huge] (0,0) node [minimum width=1cm, minimum height=1cm, above right] {$+$}; 
\end{scope}

\end{tikzpicture}
}
\end{figure}

Then the sum $D_{j_{m}}+D_{j_{m+1}}$ looks like

\begin{figure}[!ht]
\centering
\scalebox{0.6}{
\begin{tikzpicture}
 
\draw (0,0) grid (2,2);
\draw (3,0) grid (5,2);

\draw[font=\huge] (0,0.5) node [minimum width=1cm, minimum height=1cm, above left] {$\cdots$};
\draw[font=\huge] (5,0.5) node [minimum width=1cm, minimum height=1cm, above right] {$\cdots$}; 
\draw[font=\huge] (2,0.5) node [minimum width=1cm, minimum height=1cm, above right] {$\cdots$}; 
\draw[font=\huge] (-1,0.5) node [minimum width=1cm, minimum height=1cm, above left] {$D_{j_{m}}+D_{j_{m+1}}$}; 
\draw[font=\huge] (1,0) node [minimum width=1cm, minimum height=1cm, above right] {$+$}; 
\draw[font=\huge] (1,1) node [minimum width=1cm, minimum height=1cm, above right] {$+$}; 
\draw[font=\huge] (3,0) node [minimum width=1cm, minimum height=1cm, above right] {$+$}; 
\draw[font=\huge] (3,1) node [minimum width=1cm, minimum height=1cm, above right] {$+$}; 
\draw[font=\huge] (4,1) node [minimum width=1cm, minimum height=1cm, above right] {$+$};
\draw[font=\huge] (0,0) node [minimum width=1cm, minimum height=1cm, above right] {$+$};

\end{tikzpicture}
}
\end{figure}

It's easy to see that $D_{j_{m}}+D_{j_{m+1}}$ is obtained from $D_{j_{m}}$ by applying the
inverse chute move-2. Consequently, $D_{j_{m}}+D_{j_{m+1}}$ represent the same permutation. Thus,
the partition $\mathcal{P}(w)=\mathcal{P}_{s}(w)\sqcup\mathcal{P}_{I}(w)\sqcup\mathcal{P}_{\varnothing}(w)$,
constructed for $mitosis_{i}$, is also preserved by $mitosis'_{i}$.

Now we are ready to formulate the main theorem of this paper:

$\textbf{Theorem 4.1:}$ If the condition $l(ws_{i})=l(w)-1$ is sitisfied, then the following equalities take place:

1)  $\bigsqcup\limits_{D\in\mathcal{P}_{s}(w)}mitosis'_{i}(D)=\mathcal{P}(ws_{i})$

2) $\bigsqcup\limits_{D\in\mathcal{P}_{I}(w)}mitosis'_{i}(D)=\mathcal{P}_{\varnothing}(w)$.

$\textbf{Proof:}$ The disjunctivity of the unions on the left sides of the equalities 1) and 2)
is obvious: indeed, each element of the image of $mitosis'_{i}$ agrees with its preimage everywhere,
exept $i$-th and $(i+1)$-th rows. In these two adjacent rows they also coincide to the right
from the leftmost column of the type $\plel$ or $\elel$. The rest of the diagram is restored uniquely,
according to the corresponding algorithm.

Remind that if condition $l(ws_{i})=l(w)-1$ is satisfied, then $\mathfrak{G}_{ws_{i}}(x)$ is obtained from $\mathfrak{G}_{w}(x)$
by applying an the operator $\pi_{i}$.

Consider now arbitrary element $D$ of the $\mathcal{P}_{s}(w)$, such that $J=|\mathcal{J}_{i}(D)|$.
Then $x^{D}=x^{\mathcal{J}}x^{D'}$, where $D'$ is a pipe dream, obtained from $D$ by removing
all crosses in the columns with indices, belonging to $\mathcal{J}_{i}(D)$. We also have the following
equality:

$$\sum\limits_{E\in mitosis'_{i}(D)}(-1)^{|E|-l(ws_{i})}x_{E}=
(-1)^{|D|-l(w)}\left(\sum\limits_{d=1}^{J}x_{i}^{J-d}x_{i+1}^{d-1}-\sum\limits_{d=1}^{J-1}x_{i}^{J-d}x_{i+1}^{d}\right)x_{D'}=$$
$$=(-1)^{|D|-l(w)}\pi_{i}(x_{i}^{J})x_{D'}.$$

Now, if $\tau_{i}(D)=D$, then $x_{D'}$ is symmetric in the variables $x_{i}$ and $x_{i+1}$. Consequently,

$$(-1)^{|D|-l(w)}\pi_{i}(x_{i}^{J})x_{D'}=(-1)^{|D|-l(w)}\pi_{i}(x_{i}^{J}x_{D'})=\pi_{i}((-1)^{|D|-l(w)}x_{D}).$$

If $\tau_{i}(D)\neq D$, then $x_{D'}+s_{i}(x_{D'})$ is symmetric in the variables $x_{i}$ and $x_{i+1}$
and for this sum we have

$$(-1)^{|D|-l(w)}\pi_{i}(x_{i}^{J})(x_{D'}+s_{i}(x_{D'}))1=(-1)^{|D|-l(w)}\pi_{i}(x_{i}^{J}(x_{D'}+s_{i}(x_{D'})))=$$
$$=\pi_{i}((-1)^{|D|-l(w)}(x_{D}+x_{\tau_{i}(D)})).$$

(Recall, that $\tau_{i}(D)$ also lies in $\mathcal{P}_{s}(w)$ and $|D|=|\tau_{i}(D)|$)

Thus, by grouping the elements of $D\in\mathcal{P}_{s}(w)$ in accordance with the involution $\tau_{i}$,
we obtain the following equality:
$$$$
$$\sum\limits_{E\in mitosis'_{i}(\mathcal{P}_{s}(w))}(-1)^{|E|-l(ws_{i})}x_{E}=\pi_{i}(\sum\limits_{D\in\mathcal{P}_{s}(w)}(-1)^{|D|-l(w)}x_{D}).\eqno(4.1)$$
$$$$
The analogous equality can also be obtained in case of $\mathcal{P}_{I}(w)$. Now if $D\in\mathcal{P}_{\varnothing}(w)$,
i.e. $D'=D$ and $\tau_{i}(D)=D$, then $D$ is symmetric in the variables $x_{i}$ and $x_{i+1}$ and, consequently,

$$(-1)^{|D|-l(w)}\pi(x_{D})=(-1)^{|D|-l(w)}x_{D}.$$

Also if $\tau_{i}(D)\neq D$ then the sum $x_{D}+x_{\tau_{i}(D)}$ is symmetric in the variables
$x_{i}$ and $x_{i+1}$ and, consequently,

$$(-1)^{|D|-l(w)}\pi(x_{D}+x_{\tau_{i}(D)})=(-1)^{|D|-l(w)}(x_{D}+x_{\tau_{i}(D)}).$$

By groupping the elements of $\mathcal{P}_{\varnothing}(w)$ with accordance with the involution $\tau_{i}$,
we obtain the following equality
$$$$
$$\pi_{i}\left(\sum\limits_{D\in\mathcal{P}_{\varnothing}(w)}(-1)^{|D|-l(w)}x_{D}\right)=\sum\limits_{D\in\mathcal{P}_{\varnothing}(w)}(-1)^{|D|-l(w)}x_{D}.\eqno(4.2)$$
$$$$

Using these three equalities, we obtain the following result:

$$$$
$$\sum\limits_{E\in mitosis'_{i}(\mathcal{P}_{s}(w))}(-1)^{|E|-l(ws_{i})}x_{E}+\sum\limits_{E\in mitosis'_{i}(\mathcal{P}_{I}(w))}(-1)^{|E|-l(ws_{i})}x_{E}+
\sum\limits_{D\in\mathcal{P}_{\varnothing}(w)}(-1)^{|D|-l(w)}x_{D}=\eqno(4.3)$$
$$=\pi_{i}(\mathfrak{G}_{w}(x))=\mathfrak{G}_{ws_{i}}(x)=\sum\limits_{B\in\mathcal{P}(ws_{i})}(-1)^{|B|-l(ws_{i})}x_{B}.$$
$$$$

Now in order to prove 2) we only need to construct 
for an arbitrary $D$ from $\mathcal{P}_{\varnothing}(w)$ a preimage in $\mathcal{P}_{I}(w)$. Consider
the first column in the corresponding adjacent pair of rows, which is different from $\plpl$ and $\elpl$. 
If its a $\plel$ column, then, by applying
the transformation, inverse to the one on the page $11$, we will transform each $\elpl$ column on the left to the type $\plel$ and
thereby obtain the required element of $\mathcal{P}_{I}(w)$. And if it's $\elel$, then,
by applying the inverse chute move-1, we will transform $D$ to $B\backslash(i,j_{min})$, where $B$ is a pipe dream with nonempty
set $mitosis'_{i}(B)$. It follows from the theorem $\textbf{2.1}$ that $B$ belongs to $\mathcal{P}_{I}(w)$. 
Obviously, $B$ is a preimage of $D$.
Thus the action of $mitosis'_{i}$ on $\mathcal{P}_{I}(w)$ is surjective so 2) is proved. This means, that the second and third
sums on the left side of the equality $(4.3)$ are canceling out. By inputting the value $(-1)$ to the equality 
we will obtain the proof of 1). $\blacksquare$

Thus, the constructed algorithm allows us to obtain all elements of $\mathcal{P}(ws_{i})$ from the set $\mathcal{P}(w)$
in case $l(ws_{i})=l(w)-1$, giving us the inductive method of constructing Grothendieck polynomials, but now in terms of pipe dreams.

$\textbf{Acknowledgements.}$ The author is grateful to his scientific advisor Evgeny Smirnov for useful discussions.

$\addcontentsline{toc}{section}{Bibliography}$

%%%%%%%%%%%%%%%%%%%%%%%%%%%%%%%%%%%%%%%%%%%%%%%%%%%%%%%%%%%%%%%%%%%%%%%%%
\end{document}